\documentclass[doc,12pt]{apa}

\usepackage[dvips]{graphicx} 

\usepackage{amssymb, latexsym, amsmath, amscd, array, 
}

\newcommand\N {{\mathbb N}} 

\newcommand\R {{\mathbb R}} 

\newcommand\st {{\rm st}}

\title{Zooming in on infinitesimal $1-.9..$ in a post-triumvirate era}

\author{Karin Usadi Katz and Mikhail G. Katz}

\affiliation{Bar Ilan University, Israel}

\abstract{The view of infinity as a metaphor, a basic premise of
modern cognitive theory of embodied knowledge, suggests in particular
that there may be alternative ways in which one could formalize
mathematical ideas about infinity.  We discuss the key ideas about
infinitesimals via a proceptual analysis of the meaning of the
ellipsis~``$\ldots$'' in the real formula \hbox{$.999\ldots = 1$}.
Infinitesimal-enriched number systems accomodate quantities in the
half-open interval~$[0,1)$ whose extended decimal expansion starts
with an unlimited number of repeated digits $9$.  Do such quantities
pose a challenge to the unital evaluation of the symbol
``$.999\ldots$''?  We present some non-standard thoughts on the
ambiguity of the ellipsis, in the context of the cognitive concept of
{\em generic limit\/} of B.~Cornu and D.~Tall.  We analyze the
vigorous debates among mathematicians concerning the idea of
infinitesimals.

}

\rightheader{Zooming in on infinitesimal $1-.9..$ in a post-triumvirate era}

\leftheader{Karin Katz and Mikhail Katz}

\begin{document}

\maketitle


\section{1. Introduction}

Prior to the creation of operative infinitesimal-enriched continua in
the 20th century,
%
%
many mathematicians thought of infinitesimals in terms of ``naive
befogging" and ``vague mystical ideas".  Thus, R.~Courant (1937,
p.~101) presents a rather dim perspective on infinitesimals.  On that
single page, Courant manages to describe infinitesimals as
\begin{itemize}
\item
incompatible with the clarity of ideas; 
\item
entirely meaningless; 
\item
vague mystical ideas; 
\item
fog which hang round the foundations;
\item
hazy idea.
\end{itemize}

Historically speaking, infinitesimals have been difficult to
conceptualize, while a metaphor conveying the idea of proximity, such
as {\em approaching or getting closer\/}, has often been used to make
sense of their elusive properties.  In this spirit, Roh 
%
%
notes that
\begin{quote}
the words used in expressing the dynamic imagery of limits, such as
`approaching' or `getting close to', do not precisely convey the
mathematical meaning of the concept of limit \hbox{. . .}  These
expressions instead convey an everyday sense of such words, without
the precision of the mathematically rigorous meaning of limits (Roh,
2008, p.~218).
\end{quote}

We would have been more comfortable with Roh's conclusion had Roh
written, in place of the phrase {\em mathematically rigorous meaning
of limit\/}, the more accurate phrase {\em Weierstrassian definition
of limit\/}, given the existence of an alternative, {\em dynamic\/},
definition in an infinitesimal-enriched number system.%
\footnote{See footnote~\ref{cluster} for a discussion of such an
approach in the context of an infinitesimal-enriched continuum.}
R.~Ely (2010) has argued that student nonstandard conceptions,%
\footnote{Such as the conception of an infinitesimal as represented by
``$.000\ldots 1$'' (with infinitely many zeros preceding the digit
$1$) being the difference between $1$ and ``$.999\ldots 9$'' (with
infinitely many digits $9$); cf.~Lightstone's notation in
footnote~\ref{nine}.}
routinely perceived and dismissed as erroneous conceptions, can
actually serve as valuable tools in learning calculus.

We propose an approach to infinitesimals via an examination of the
meaning of the real formula $.999\ldots=1$, henceforth referred to as
the {\em unital evaluation\/} of the symbol~``$.999\ldots$''.

In Section~2, we examine the cognitive concept of {\em generic
limit\/} of Cornu and Tall, in relation to a hyperreal approach to
limits, and exhibit a hyperreal quantity in~$[0,1)$ with an unlimited
number of repeated~$9$s.  In Section~3, we represent such a quantity
graphically by means of an infinite-magnification microscope, already
exploited for pedagogical purposes by Keisler (1986)
%
%
and Tall (1980),
%
%
and use it to calculate~$f'(1)$ where~$f(x)=x^2$, in the framework of
Tall's locally straight approach.  The historical Section~4 contains
an examination of the views of R.~Courant, I.~Lakatos, E.~Bishop, and
A.~Heyting as pertaining to infinitesimals.  In Section~5, we develop
an applied-mathematical model of a hypercalculator so as to explain a
familiar phenomenon of a calculator returning a string of 9s in place
of an integer.  In Section~6 we propose that research be undertaken to
see if infinitesimals are a better support for students to come to
grips with the notion of limits, and draw some conclusions regarding
the foundations of analysis and the ramifications of the use of an
infinitesimal-enriched continuum in teaching.%
\footnote{See also Sad, Teixeira, Baldino (2001, p.~286).}

\section{2. Generic limit}
\label{two}

D. Tall (2009)
%
%
notes that ``the infinite decimal~$0.999\ldots$ is intended to
signal the {\em limit\/} of the sequence~$0.9,\;$~$0.99,\;$~$0.999$,
... which is~$1$, but in practice it is often imagined as a limiting
process which never quite reaches~$1$.''
Tall (1991, 2009)
%
%
describes a concept in cognitive theory he calls a {\em generic limit
concept\/} in the following terms: 
\begin{quote}
[I]f a quantity repeatedly gets smaller and smaller and smaller
without ever being zero, then the limiting object is naturally
conceptualised as an extremely small quantity that is not zero (Cornu,
1991).
%
%
Infinitesimal quantities are natural products of the human imagination
derived through combining potentially infinite repetition and the
recognition of its repeating properties.
\end{quote}

Nonstandard student conceptions about infinitesimals were recently
analyzed by R.~Ely (2007, 2010).
%
%
J.~Monaghan (2001, p.~248),
%
%
based on field studies, concluded as follows: ``do infinite numbers of
any form exist for young people without formal mathematical training
in the properties of infinite numbers?  The answer is a qualified
`yes'.''  Numbers with infinitely many digits are not beyond the
intellectual capacity of children.

An examination of the meaning of the real formula $.999\ldots=1$ in
the context of an infinitesimal-enriched number system, suggests a
mathematical realisation of the cognitive concept of {\em generic
limit\/}, in terms of a choice of a hypernatural number, thought of as
an encapsulation in the context of a metaphor of an ``ever larger
natural number'', as follows.  In the familiar finite domain, we
evaluate the formula
\begin{equation*}
1+r+r^2+\ldots+r^{n-1}=\frac{1-r^{n}}{1-r}
\end{equation*}
at~$r=\frac{1}{10}$, we obtain
\begin{equation*}
1+\frac{1}{10}+\frac{1}{100}+\ldots + \frac{1}{10^{n-1}}=
\frac{1-\frac{1}{10^{n}}}{1-\frac{1}{10}} ,
\end{equation*}
or alternatively

\begin{equation*}
\underset{n}{\underbrace{1.11\ldots 1}} =
\frac{1-\frac{1}{10^{n}}}{1-\frac{1}{10}} .
\end{equation*}
Multiplying by~$\frac{9}{10}$, we obtain
\begin{equation*}
\begin{aligned}
.\underset{n}{\underbrace{999\ldots 9}}\; & = \frac{9}{10} \left(
\frac{1-\frac{1}{10^{n}}}{1-\frac{1}{10}} \right) \\&=1 -
\frac{1}{10^{n}} 
\end{aligned}
\end{equation*}
for every~$n\in \N$.  As~$n$ increases without bound, the formula
\begin{equation}
\label{21}
.\underset{n}{\underbrace{999\ldots 9}}\; = 1- \frac{1}{10^{n}}
\end{equation}
becomes
\begin{equation*}
.999\ldots = 1;
\end{equation*}
at any rate, so goes the traditional mathematical account of the
matter of~``$.999\ldots$''.

Cognitively speaking, the ``underbrace'' formula~\eqref{21} suggests
exploiting a generic (in Tall's sense) ``infinite natural number'',
which we will denote $[\N]$, encapsulating, as it were, the ever
increasing sequence 
\[
\langle \ldots, n-1, \; n, \; n+1, \ldots \rangle,
\]
so as to obtain
\begin{equation}
\label{22}
.\underset{[\N]}{\underbrace{999\ldots }}\; = 1-
\frac{1}{10_{\phantom{I}}^{[\N]}}.
\end{equation}
The latter formula is suggestive of an infinitesimal difference in the
minute amount of~$\frac{1}{10_{\phantom{I}}^{[\N]}}$, with the
quantity $.\underset{[\N]}{\underbrace{999\ldots }}$ falling just
short of~$1$.  The procedure presented so far can be described as
inhabiting the {\em proceptual\/%
\footnote{The term {\em procept\/} was coined by Gray and Tall
(1991).}
%
symbolism\/} world, namely the second of Tall's three worlds (see
Tall, 2010).
%
%

A remarkable passage by Leibniz is a testimony to the enduring appeal
of the metaphor of infinity, even in its, paradoxically, {\em
terminated\/} form.  In a letter to Johann Bernoulli dating from june
1698, Leibniz speculated concerning ``lines . . . which are terminated
at either end, but which nevertheless are to our ordinary lines, as an
infinite to a finite'' (see Jesseph, 1998).
%
%
And again, he speculates as to the possibility of
\begin{quote}
a point in space which can not be reached in an assignable time by
uniform motion.  And it will similarly be required to conceive a time
terminated on both sides, which nevertheless is infinite, and even
that there can be given a certain kind of eternity . . . which is
terminated.
\end{quote}
To be sure, Leibniz rejected such metaphysics, and ultimately
conceived of both infinitesimals and infinite quantities as ideal
numbers, falling short of a metaphysical reality of familiar
appreciable quantities.  To formalize such quantities mathematically,
the transition to Tall's third world of {\em axiomatic formalism\/}
can be accomplished through a ``generic'' analysis of the classical
construction of the reals as equivalence classes of Cauchy sequences
of rational numbers, that parallels our ``generic'' analysis of the
real formula $.999\ldots=1$.  Here one first {\em relaxes\/} the
equivalence relation, so that certain null sequences (i.e.~sequences
tending to zero) will become inequivalent (and represent distinct
infinitesimals).  Furthermore, one {\em extends\/} the relation to
arbitrary (not merely Cauchy) sequences of rational numbers.  Thus,
the sequence
\[
\langle u_n, \; n=1, 2, 3, \ldots \rangle = \langle \N \rangle,
\]
listing all the natural numbers in increasing order, will be (at the
level of the new equivalence class) encapsulated as an infinite
hypernatural number, denoted $[\N]$.%
\footnote{Such a construction is known as the ultrapower construction
and is due to Luxemburg in 1962, see (Luxemburg 1964) and (Goldblatt,
1998).  To obtain the full hyperreal field, one starts with sequences
of real numbers.}
%

The resulting alternative evaluation \eqref{22} of the symbol
``$.999\ldots$'' thus employs the string of natural numbers $\langle
1,2,3,\ldots \rangle$, and could be accordingly termed its {\em
natural string\/} evaluation, a ``generic'' competitor to unital
evaluation.%
\footnote{\label{cluster}Each real number is accompanied by a cluster
(alternative terms are prevalent in the literature, such as {\em
halo\/}, but the term {\em cluster\/} has the advantage of being
self-explanatory) of hyperreals infinitely close to it.  The standard
part function collapses the entire cluster back to the standard real
contained in it.  The cluster of the real number~$0$ consists
precisely of all the infinitesimals.  Every infinite hyperreal
decomposes as a triple sum $H+r+\epsilon$, where~$H$ is a hyperinteger
(see below),~$r$ is a real number in~$[0,1)$, and $\epsilon$ is
infinitesimal.  Varying~$\epsilon$ over all infinitesimals, one
obtains the cluster of~$H+r$.  A hyperreal number~$H$ equal to its own
integer part: $H = [H]$ is called a hyperinteger.  Here the integer
part function is the natural extension of the real one.  Limited
(finite) hyperintegers are precisely the standard ones, whereas the
unlimited (infinite) hyperintegers are sometimes called non-standard
integers.  The {\em limit\/} $L$ of a convergent sequence $\langle u_n
\rangle$ is the standard part st of the value of the sequence at an
infinite hypernatural: $L=\st(u_H)$, for instance at~$H=[\N]$.  In
connection with $.999\ldots$, I.~Stewart notes that 

\begin{quote}
[t]he standard analysis answer is to take `$\ldots$' as indicating
passage to a limit.  But in non-standard analysis there are many
different interpretations (Stewart, 2009, p.~176).
\end{quote}

Some additional details to be found in (Katz \& Katz, 2010).
Historically, there have been two main approaches to infinitesimals.
The approach of Leibniz postulates the existence of infinitesimals of
arbitrary order, while B.~Nieuwentijdt favored nilpotent (nilsquare)
infinitesimals (see J.~Bell 2009).  Bell notes that a Leibniz
infinitesimal is implemented in the hyperreal continuum of Robinson,
whereas a Nieuwentijdt infinitesimal is implemented in the smooth
infinitesimal analysis of F.W.~Lawvere, based on intuitionistic
logic.}
Ultimately, the traditional insistence on the unital evaluation of the
string with an unlimited number of~$9$s, is necessitated by the
absence of infinitesimals in the traditional number system.

\section{3. Infinitesimals under magnifying glass}

Tall's ``locally straight'' approach employing computer graphics and
the metaphor of an infinite microscope, can also be used effectively
in the context of an analysis of the real equality $.999\ldots=1$, so
to give students a geometric feel for the process which is occurring
just short of $1$.  Recall that the symbol ``$\infty$'' is employed in
standard real analysis to define a formal completion of the real
line~$\R$, namely
\begin{equation}
\label{51}
\R\cup \{\infty\}.
\end{equation}

\begin{figure}
\includegraphics[height=4.8in]{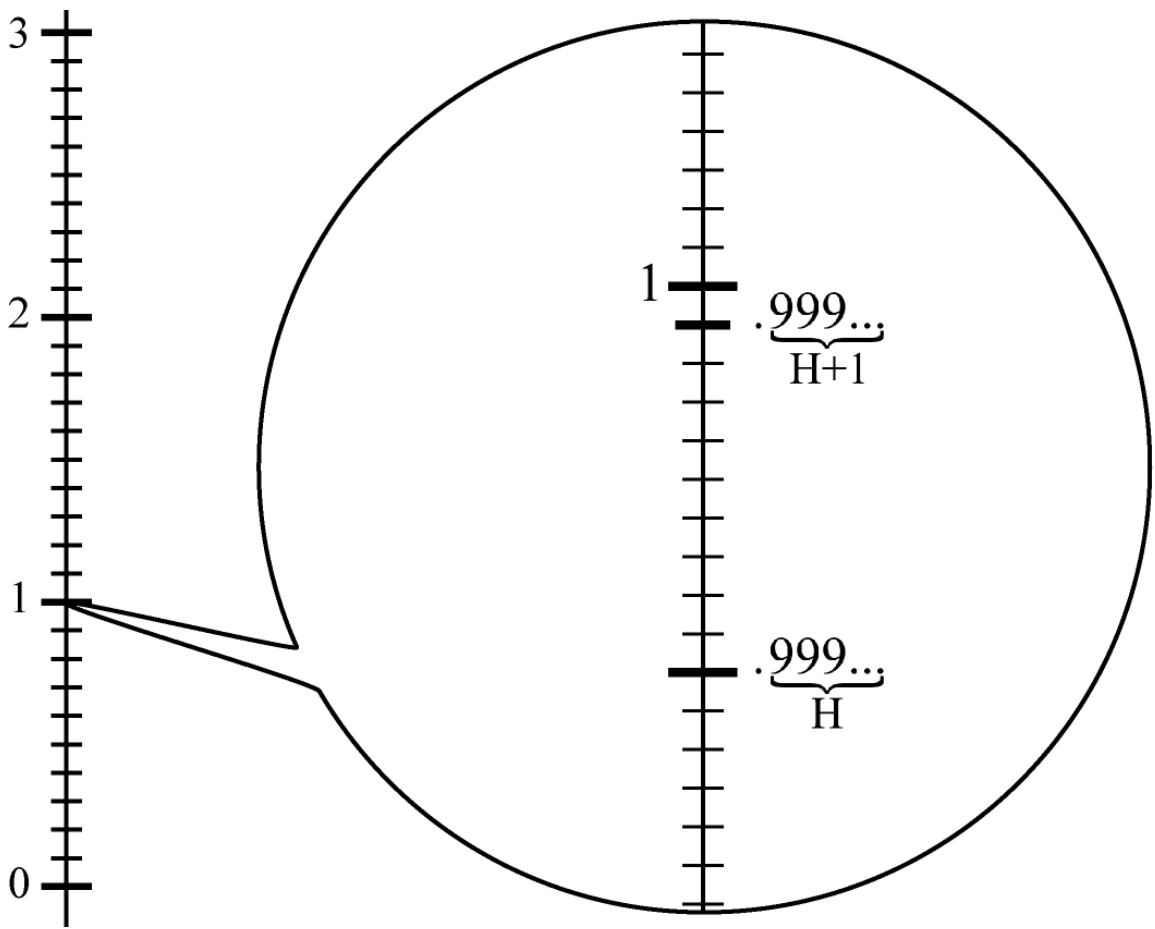}
\caption{Three infinitely close hyperreals under a microscope}
\label{micro}
\end{figure}

Sometimes a formal point~``$-\infty$'' is added, as well.  Such a
formal device is helpful in simplifying the statements of certain
theorems (which would otherwise have a number of subcases).%
\footnote{\label{ftwo}The symbol is used in a {\em different\/} sense
in projective geometry, where adding a point at infinity~$\{\infty\}$
to~$\R$ results a circle: $\R\cup \{\infty\} \approx S^1$.}

We have refrained from exploiting the symbol ``$\infty$'' to denote an
infinite hypernatural number, as in
\[
.\underset{\infty}{\underbrace{999\ldots}},
\]
so as to avoid the risk of creating a false impression of the
uniqueness, as in \eqref{51}, of such an infinite point (see
Figure~\ref{micro}).  We represent the quantity
\begin{equation}
\label{52}
.\underset{H}{\underbrace{999\ldots}}
\end{equation}
visually by means of an infinite-resolution microscope already
exploited for pedagogical purposes by Keisler (1986).
%
%
Similarly, we can follow Tall's ``locally straight'' approach, and
exploit our ``generic'' analysis of ``$.999\ldots$'' in order to
calculate the slope of the tangent line to the curve~$y=x^2$ at $x=1$.

\begin{figure}
\includegraphics[height=2.498in]{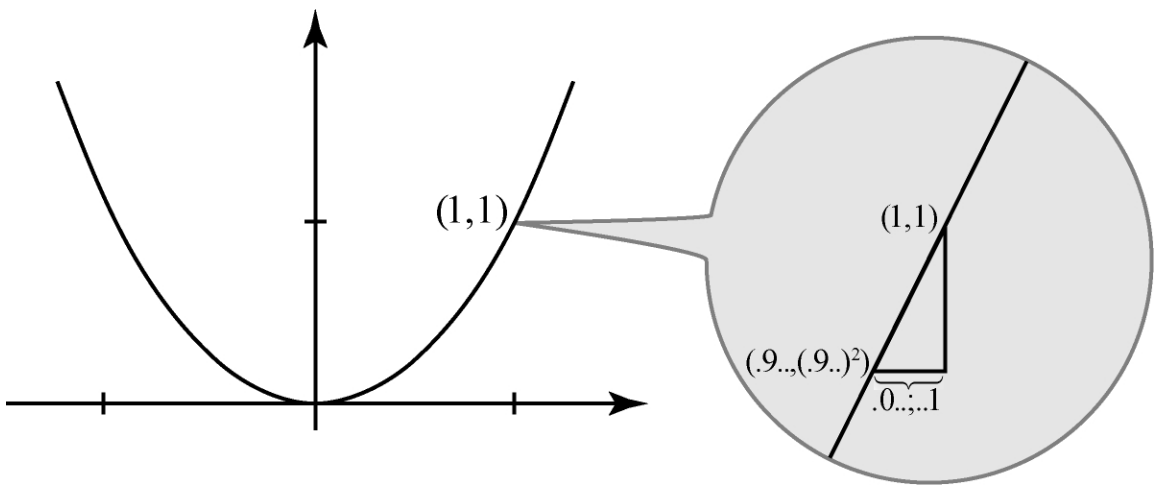}
\caption{Enlargement of an infinitesimal segment of a parabola as
calculation of the slope~$f'(1)$ for $f(x)=x^2$.}
\label{slope}
\end{figure}

Thus, we first compute the ratio
\[
\frac{\Delta y}{\Delta x} = \frac{(.9..)^2 -1^2}{.9.. - 1} =
\frac{(.9.. -1)(.9.. + 1)}{.9.. -1} = .9.. + 1,
\]
where we have deleted the underbrace
$\underset{H}{\underbrace{\phantom{.999\ldots}}}$ and also shortened
the symbol~$.999\ldots$ to~$.9..$, so as to lighten the notation for
the purposes of this section only.  Next, we note that the ``standard
part'' function, denoted ``st'', associates to each finite quantity,
the unique real number infinitely close to it.  The slope is then
computed as the standard part%
\footnote{\label{nine}Note that $\Delta x = .9.. - 1 = -.0..;..01$ in
Lightstone's notation (see Lightstone, 1972),
%
%
where the digit ``1'' appears at infinite decimal rank~$H$.}
\[
\frac{dy}{dx} = \st\left( \frac{\Delta y}{\Delta x} \right) =
\st(.9.. + 1) = \st(.9..) + 1 = 1+1=2,
\]
as illustrated in Figure~\ref{slope}.

In the context of geometric visualisation of infinitesimals, perhaps a
helpful parallel is provided by the famous animated film {\em
Flatland\/}, based on the classic book (Abbott, 2008).
%
%
Here the two--dimensional creatures are unable to conceive of what we
think of as the sphere in 3-space, due to their dimension limitation.
Similarly, one can conceive of the difficulty in the understanding of
the unital evaluation of~``$.999\ldots$'', as due to the limitation of
the standard real vision.

The flatlanders see the sphere, as it cuts across their existential
plane, as the creation of a point, continuously growing through a
family of concentric circles, and then decreasing in radius and
petering out as the circle shrinks to a point and disappears.
Meanwhile, E. A. Abbott sees, of course, a sphere as it moves across a
plane in~$3$-space.

Over half a millenium earlier, in 1310, Yitzchak HaYisraeli (1310)
%
%
may have seen more than what Mr.~Abbott saw.  Before the disappearance
of the shrinking circle, he perceives ``the smallest of them, a cap
very, very small so you would think there is no cap smaller than it,
with the point of the pole fixed, as it were, in its midst as its
center'' (HaYisraeli, 1310, p.~10).%
\footnote{See (Goldwurm, 2001, p.~104) for biographical information on
HaYisraeli.}
Note that he mentions the pole as the center of the {\em very, very
small cap}, as distinct from the cap itself.  This indicates the the
cap has not yet degenerated to a point in the author's mind.  Yet it
is {\em very, very small\/} and {\em there is no cap smaller than it},
indicating that he is referring to an infinitesimal cap.

Perhaps Yitzchak HaYisraeli can be thought of as a precursor of the
{\em generic limit\/}, all this a century prior to Nicholas of Cusa.
Nicholas of Cusa (15th century) considered the circle as a polygon
with an infinite number of sides, inspiring Kepler to formulate his
{\em bridge of continuity\/} (see Tall, 2010).
%
%
Yitzchak HaYisraeli's book {\em Yesod Olam\/} was written over a
decade before the birth of Nicole Oresme, another pioneer of
infinitesimals.  HaYisraeli (1310, p.~12 and fol.)  also presents a
detailed statement of the spherical law of sines and both cases of the
spherical law of cosines for right-angled triangles,
%
%
and mentions, as one of his illustrious predecessors at Toledo, the
astronomer Ibrahim al-Zarqali (see HaYisraeli, 1310).
%
%

That the notion of an infinitesimal appeals to intuition is
acknowledged even by its critics, see R.~Courant's comment in
Section~4.
The infinitesimal would not go away inspite of what is, by now, over a
century of~$(\epsilon, \delta)$-ideology, see Bishop's comment in
Section~4.

\section{4. Second thoughts about infinitesimals}

In the introduction, we quoted Courant's unflattering appraisal
(predating Robinson's work) of infinitesimals.  Yet, Courant
acknowledges the appeal of the infinitesimal: ``For a great many
simple-minded people it undoubtedly has a certain charm, the charm of
mystery which is always associated with the word `infinite'.''
Courant was unable to peer through the hazy mystical fog the way
Robinson would.  It should be kept firmly in mind that Courant's
criticism {\em predated\/} Robinson's theory, unlike more recent
criticisms.  Courant's criticism was not without merit {\em at the
time of its writing\/}, namely, a quarter century prior to Robinson's
work.  I.~Lakatos
%
%
wrote in 1966 as follows:
\begin{quote}
Robinson's work . . . offers a rational reconstruction of the
discredited infinitesimal theory which satisfies modern requirements
of rigour and which is no weaker than Weierstrass's theory.  This
reconstruction makes infinitesimal theory an almost respectable
ancestor of a fully fledged, powerful modern theory, lifts it from the
status of pre-scientific gibberish, and renews interest in its partly
forgotten, partly falsified history (Lakatos, 1978, p.~44).
\footnote{\label{ehrlich}The falsification problem is analyzed in
Ehrlich (2006).}
%
\end{quote}

Not everyone was persuaded by Lakatos' appreciation of the
significance of Robinson's work.  A decade later, in 1975, Courant's
{\em duty of avoiding every such hazy idea\/} was taken up under a
constructivist banner by E.~Bishop, affectionately described as the
``reluctant guru'' of constructive mathematics by his advisor,
P.~Halmos (1985, p.~162).
%
%
In his essay (Bishop, 1975)
%
%
cast in the form of an imaginary dialog between Brouwer and Hilbert,
Bishop anchors his foundational stance in a species of mathematical
constructivism.  Thus, Bishop's opposition to Robinson's
infinitesimals, expressed in a vitriolic%
\footnote{\label{vit}Historians of mathematics have noted the
vitriolic tenor of Bishop's criticism, see e.g.~(Dauben, 1996),
p.~139.}
%
review  (Bishop, 1977)
%
%
of Keisler's textbook,%
\footnote{Namely, review of an earlier edition of (Keisler, 1986).}
%
was to be expected (and in fact was anticipated by editor Halmos).

Two years earlier, Bishop had already expressed his views about
non-standard analysis and its use in teaching in a brief paragraph
toward the end of his essay ``Crisis in contemporary mathematics'',
see Bishop (1975, p.~513-514).
%
%
Having discussed Hilbert's formalist program, he writes:
\begin{quote}
A more recent attempt at mathematics by formal finesse%
\footnote{The description of Hilbert's program as ``formal finesse''
has been objected to by many authors.  Avigad and Reck (2001)
%
%
provide a detailed discussion of the significance, and {\em
meaning\/}, of Hilbert's program.}
is non-standard analysis. I gather that it has met with some degree of
success, whether at the expense of giving significantly less
meaningful proofs I do not know. My interest in non-standard analysis
is that attempts are being made to introduce it into calculus
courses. 
\end{quote}
Bishop concludes: ``It is difficult to believe that debasement of
meaning could be carried so far.''  Bishop's view of the introduction
of non-standard analysis in the classroom as no less than a {\em
debasement of meaning\/}, was duly noted by Dauben~(1992).
%
%

Bishop's sentiments toward non-standard calculus stand in sharp
contrast with those of his fellow intuitionist A. Heyting (1973,
p.~136)
%
%
who felt that
\begin{quote}
[Robinson] connected [an] extremely abstract part of model theory with
a theory apparently so far apart as the elementary calculus.  In doing
so [he] threw new light on the history of the calculus by giving a
clear sense to Leibniz's notion of infinitesimals.%
\footnote{\label{f49}Note that Schubring (2005)
attributes the first systematic use of infinitesimals as a
foundational concept, to Johann Bernoulli.}
\end{quote}

Bishop clarified what exactly it was he had in mind when he spoke of
{\em debasement of meaning\/} in an earlier text (Bishop, 1985)
%
%
distributed in 1973 and eventually published in 1985.  Here Bishop
writes: ``Brouwer's criticisms of classical mathematics were concerned
with what I shall refer to as `the debasement of meaning' '' (1985,
p.~1).
%
%
In Bishop's own words, the {\em debasement of meaning\/} expression,
employed in his {\em Crisis\/} text to refer to non-standard calculus,
is in fact a criticism of {\em classical mathematics\/} as a whole,
exposing a thinly disguised foundational agenda in his criticism of
non-standard calculus as well.

\label{fef}
To illustrate how Bishop anchors his foundational stance in a species
of mathematical constructivism, note that
he writes:
\begin{quote}
To my mind, it is a major defect of our profession that we refuse to
distinguish . . . between integers that are computable and those that
are not . . . the distinction between computable and non-computable,
or constructive and non-constructive is the source of the most famous
disputes in the philosophy of mathematics...  (Bishop, 1975,
pp.~507-508).
%
%
\end{quote}

On page 511, Bishop defines a principle (LPO) as the statement that
``it is possible to search `a sequence of integers to see whether they
all vanish' '', and goes on to characterize the dependence on the LPO
as a procedure both Brouwer and Bishop himself reject. S.~Feferman
%
%
explains the principle as follows:
\begin{quote}
Bishop criticized both non-constructive classical mathematics and
intuitionism.  He called non-constructive mathematics ``a scandal",
particularly because of its ``deficiency in numerical meaning".  What
he simply meant was that if you say something exists you ought to be
able to produce it, and if you say there is a function which does
something on the natural numbers then you ought to be able to produce
a machine which calculates it out at each number (Feferman, 2000).
\end{quote}
Elsewhere, Feferman identifies LPO as a special case of the law of
excluded middle.

It recently came to light (Manning, 2009)
%
%
that Bishop never uttered his criticism of non-standard calculus in
his oral presentation at the 1974 workshop (which helps explain the
absence of any critical reaction to such {\em debasement\/} on the
part of the audience in the discussion session, included at the end of
the published version of his talk), but rather inserted it at the
galley proof stage of publication.  Bishop fails to acknowledge in his
review in the {\em Bulletin\/} that his criticism is motivated by his
foundational preoccupation with the law of excluded middle, and with
what he calls {\em numerical meaning\/} in his {\em Crisis\/} essay.
In short, Bishop is criticizing apples for not being oranges: the
critic (Bishop) and the criticized (non-standard analysis) do not
share a common foundational framework.  Note that a similar point was
mentioned by M.~Davis (1977, p.~1008).
%
%
Note that the foundational framework of non-standard analysis, namely
the Zermelo-Frankel set theory with the axiom of choice (ZFC), is the
framework of the vast majority of the readers of the {\em Bulletin},
at variance with Bishop's preferred framework, in which his {\em
debasement of meaning\/} criticism would apply equally well to all of
mainstream mathematics.  This point was alluded to in Feferman's
comment cited above.

What Bishop sees as {\em debasement of meaning\/} in classical
mathematics is the alleged absence of {\em numerical meaning\/}
already alluded to above.  Classical logic relies on the {\em law of
excluded middle\/} (LEM), the key ingredient in proofs by
contradiction.  Bishop tends to conflate his narrow notion of {\em
numerical meaning\/} (narrowly defined as the avoidance of LEM) with
{\em meaning\/} in a wider epistemological sense.

Meanwhile, infinitesimal calculus as developed by the founders of the
discipline%
\footnote{See footnote~\ref{f49} for a historical clarification.}
possesses numerical meaning of the post-LEM variety, as expressed in
explicit formulas for the derivatives and integrals of standard
functions.  Infinitesimal calculus remains in the post-LEM category
even after Robinson's work, due to the nature of the construction of
the hyperreal number system.%
\footnote{An alternative infinitesimal-enriched intuitionistic
continuum has been developed by Lawvere, see footnote~\ref{cluster}.}

Dauben (1996, p.~135)
%
%
presents a detailed analysis, in the areas both of pedagogy and of
meaning, of Bishop's criticisms, and describes them as ultimately
``unfounded".  Dauben (1996, p.~133)
%
%
identifies an intriguing point of convergence between Bishop and
Robinson, namely that the history of the calculus has been, in
Bishop's words, ``systematically distorted to support the status quo"
(see (Bishop, 1975, p.~508)).
%
%
Dauben describes as ``one of the most important achievements of
Robinson's work in non-standard analysis", ``his conclusive
demonstration of the poverty of [the] kind of historicism" that
focuses exclusively on an alleged triumph of ``Weierstrassian
epsilontics over infinitesimals in making the calculus ``rigorous" in
the course of the 19th century".%
\footnote{The issue of the falsification of the history of the
calculus was discussed by Lakatos, see main text at
footnote~\ref{ehrlich}.}
Non-standard calculus in the classroom was analyzed in the Chicago
study by K.~Sullivan (1976).
%
%
Sullivan showed that students following the non-standard calculus
course were better able to interpret the sense of the mathematical
formalism of calculus than a control group following a standard
syllabus.  Such a conclusion was also noted by
M.~Artigue (1994, p.~172).
%
%
A synthesis of such experiments was made by Bernard Hodgson (1994)
%
%
in 1992, and presented at the ICME-7 at Quebec.  G.~Schubring (2005,
p.~153)
%
%
points out that an alternative approach to calculus developed by a
German mathematician ``has been . . .  unable to win as much celebrity
and as many adherents as [non-standard analysis]''.  Leibniz historian
H.~Bos (1974, p.~13)
acknowledged that Robinson's hyperreals provide~a
\begin{quote}
preliminary explanation of why the calculus could develop on the
insecure foundation of the acceptance of infinitely small and
infinitely large quantities.
\end{quote}
F.~Medvedev (1998)
%
further points out that nonstandard analysis ``makes it possible to
answer a delicate question bound up with earlier approaches to the
history of classical analysis.  If infinitely small and infinitely
large magnitudes are [to be] regarded as inconsistent notions, how
could they [have] serve[d] as a basis for the construction of so
[magnificent] an edifice of one of the most important mathematical
disciplines?''
\vskip.5in

\section{5. Hypercalculator returns~$.999\ldots$}

Our goal here is to illustrate how a non-standard
number~``$.999\ldots$'' can be exploited in an explanation of a
familiar phenomenon from routine calculator use.  This topic is too
advanced to be presented at a highschool level, as it involves
Newton's algorithm.  The latter is an iterative procedure for finding
a sequence~$x_0, x_1, \ldots$ rapidly converging to a zero of a smooth
function, under suitable hypotheses on the derivatives, and thus
requires a familiarity with calculus.  On the other hand, it could be
presented as an ``enrichment'' topic to a class of college students,
already familiar with standard calculus, in the context of a first
exposure to the hyperreals.

Everyone who has ever held an electronic calculator is familiar with
the curious phenomenon of it sometimes returning the value
\begin{equation*}
.999999
\end{equation*}
in place of the expected~$1.000000$.  For instance, a calculator
programmed to apply Newton's method to find the zero of a function,
may return the~$.999999$ value as the unique zero of the function
\[
\log x.
\]
Developing a model to account for such a phenomenon is complicated by
the variety of the degree of precision displayed, as well as the
greater precision typically available internally than that displayed
on the LCD.  To simplify matters, we will consider an idealized model,
called a {\em hypercalculator}, of a theoretical calculator that
applies Newton's method precisely~$H$ times, where~$H$ is a fixed
infinite hypernatural, e.g.~$[\N]$, as discussed in Section~2.

It can be established that {\em if~$f$ is a concave strictly
increasing differentiable function with domain an open
interval~$(1-\epsilon, 1+\epsilon)$ and vanishing at its midpoint,
then the hypercalculator applied to~$f$ will return a hyperreal
decimal~$.999\ldots$ with an initial segment consisting of an
unlimited number of repeated~$9$'s.}

Indeed, it is well-known that Newton's algorithm converges under the
above hypotheses.  We will reproduce the main calculation of the
standard proof, emphasizing the novelty that the strict inequality can
be retained in the end, as in \eqref{a1}.

Assume for simplicity that~$f(x_0)<0$.  We have 
\begin{equation*}
x_1= x_0 + \frac{|f(x_0)|}{f'(x_0)}.
\end{equation*}
By the mean value theorem, there is a point~$c$ such that~$ x_0 < c <
1$ where~$f'(c)= \frac{|f(x_0)|}{1-x_0}$, or
\begin{equation*}
\frac{|f(x_0)|}{f'(c)}= 1-x_0 .
\end{equation*}
Since~$f$ is concave, its derivative~$f'$ is decreasing, hence
\begin{equation*}
x_1 = x_0 + \frac{|f(x_0)|}{f'(x_0)} < x_0 +1 -x_0 = 1 .
\end{equation*}
Thus~$x_1 < 1$.  Inductively, the point~$x_{n+1}= x_n +
\frac{|f(x_n)|} {f'(x_n)}$ satisfies~$x_{n}<1$ for all~$n$.  By the
transfer principle of non-standard analysis, the hyperreal~$x_{H}$
satisfies a strict inequality
\begin{equation}
\label{a1}
x_{H}<1,
\end{equation}
as well.  Hence the hypercalculator returns a value strictly smaller
than~$1$ yet infinitely close to~$1$, and therefore starts with an
unlimited number of 9s.

\section{6. Conclusions}

The utility of infinitesimals is known to transcend the ``$.999...$''
issue.  In the fall of 2008, the second-named author taught a course
in calculus using Keisler's textbook {\em Elementary Calculus\/}
(Keisler, 1986)
%
%
to a group of 25 freshmen.  The TA had to be trained as well, as the
material was new to the TA.  The students have never been so excited
about learning calculus, according to repeated reports from the TA.
Two of the students happened to be highschool teachers (they were
somewhat exceptional in an otherwise teenage class).  They said they
were so excited about the new approach that they had already started
using infinitesimals in teaching basic calculus to their 12th graders.
After the class was over, the TA paid a special visit to the
professor's office, so as to place a request that the following year,
the course should be taught using the same approach.  Furthermore, the
TA volunteered to speak to the chairman personally, so as to impress
upon him the advantages of the approach exploiting an
infinitesimal-enriched number system.  {\em Restoring\/}%
\footnote{Infinitesimals were in classroom use as late as the 1930s,
if sometimes surreptitiously, see (Roquette 2008) and (Luzin 1931).}
%
an infinitesimal-enriched number system to the classroom will enhance
the presentation of a million calculus teachers around the globe.

The difficulty of the Weierstrassian limit concept is no secret either
to professional mathematicians, or to professional educators.  One
approach to the difficulty has been to view it as intrinsic in the
subject matter.  Thus, C. Boyer (1949, p.~298)
%
%
sees the crowning of 2500 years of ``investigations leading to the
calculus'' from Pythagoras onward, in the ``satisfactory definitions
of number and the infinite'' established by ``the great triumvirate:
Weierstrass, Dedekind, and Cantor''.  To his credit, Boyer notes on
page 299 that ``[t]here is a strong temptation on the part of
professional mathematicians and scientists to seek always to ascribe
great discoveries and inventions to single individuals.''  On page
287, we learn that ``the limit concept does not involve the idea of
{\em approaching\/}, but only a static state of affairs.''  On page
298, it turns out that ``any criticism of the use of the infinite in
defining irrational number or in the limit concept is answered by
Cantor's work, which clarifies the situation.''  Boyer asserts that
such criticisms were answered by Cantor.  Or perhaps not entirely?
Boyer equips his phrase with a footnote 92, which reads as follows:
``It should, perhaps, be observed at this point that the theory of
infinite aggregates has resulted also in a number of puzzling and as
yet unresolved antimonies.''  Nonetheless, ``the real number of
Dedekind is in a sense a creation of the human mind, independent of
intuitions of space and time'' (Boyer, 1949, p.~292).
%
%

A cognitive scientist would agree emphatically that the real number is
a creation of the human mind; he would disagree that such a creation
(as all embodied knowledge) could be independent of intuitions of
space and time.  However, on page~305 we are told that ``[t]he bases
of the calculus were then defined formally in terms only of number and
infinite aggregates, with no corroboration through an appeal to the
world of experience either possible or necessary.''

One wonders how many calculus instructors may have adopted Boyer's
view concerning a lack of necessity of appeals to the world of
experience, in formally defining the bases of the calculus.  On page
308, Boyer describes mathematics as ``the symbolic logic of possible
relations'', and is content to conclude his book on page 309 with a
final flourish, describing mathematics as ``a syllogistic elaboration
of arbitrary premises''.

Boyer has thus invested the great triumvirate with the highest
ministry in a disembodied world of mathematics, stripped of both
intuition and motion, but firmly afloat upon the twin whales of a
Cantorian set-theoretic paradise of infinite aggregates, on the one
hand, and unassailable symbolic {\em logic\/}, on the other.  Let us
therefore hear out the {\em logicians\/} themselves---the best of
them.

Abraham Robinson, speaking at Bristol in 1973, confides as follows:
``I do not believe in the primacy of set theory over all other
branches of mathematics'' (Robinson, 1975, p.~48); (Robinson, 1979,
p.~563).
%
%
But what about the infinite aggregates?  ``[M]athematical theories
which, allegedly, deal with infinite totalities do not have any
detailed . . . reference'' (Robinson, 1975, p.~42; 1979, p.~557).
%
%
Two decades before Lakoff \& N\'u\~nez (2000),
Robinson describes mathematical infinities as lacking a detailed
reference, i.e., as metaphors.  But is it conceivable that the
empirical scientists should dare throw off the yoke of the great
triumvirate?  On this point, Robinson's deference to the empirical
scientist is as palpable as it is visionary:
\begin{quote}
At this point, we should also consider the possibility of future
developments in the empirical sciences which will affect the areas
with which we are concerned.  I can think of one such development,
which will surely occur in the fullness of time, although I have no
means of judging when.  This is the possibility of analyzing in detail
the neurophysiological processes in the brain which correspond to its
mathematical activity (Robinson, 1975, p.~48-49); (Robinson, 1979,
p.~563-564).
%
%
\end{quote}

An increasing sentiment in foundational and educational circles is
that there are alternative ways in which to formalize mathematical
ideas (in this case, ideas about infinity).  It is not uncommon to
hear opinions to the effect that one needs to break the dogma that
foundations in mathematics are a `solid building', or that it is not
possible to have different explanations to justify formal results.%
\footnote{In a letter to M.~Vygodski\u\i, the mathematician Luzin
questioned whether the Weierstrassian approach to the foundations of
analysis ``corresponds to what is in the depths of our consciousness
\hbox{. . .}  I cannot but see a stark contradiction between the
intuitively clear fundamental formulas of the Integral calculus and
the incomparably artificial and complex work of the `justification'
and their `proofs' '' (Luzin,~1931).
The publication of the text {\em Fundamentals of Infinitesimal
Calculus\/}, by Vygodski\u\i, in 1931, provoked sharp criticisms.
Luzin wrote his (two) letters to counterbalance such criticisms, and
took the opportunity to elaborate his own views concerning
infinitesimals.}  
Other educators feel that mathematical ideas can be, and have been,
controversial.

\section*{Acknowledgments}

We are grateful to David Ebin and David Tall for a careful reading of
an earlier version of the manuscript, and for making numerous helpful
suggestions.  We thank the editor of the article and the anonymous
referees for an exceptionally thorough analysis of the shortcomings of
the version originally submitted, resulting, through numerous
intermediate versions, in a more focused text.

\end{document}